\documentclass[MSNbibl,number,citesort,seceqn,dvips]{arxbj}
\usepackage{upgreek}
\usepackage{graphicx}

\aid{0}
\volume{18}
\issue{3}
\pubyear{2012}
\firstpage{836}
\lastpage{856}
\doi{10.3150/11-BEJ365}

\makeatletter
\newcommand{\A}{{ A}}
\newcommand{\B}{{ B}}
\newcommand{\C}{{ C}}
\newcommand{\D}{{ D}}
\newcommand{\G}{{ G}}
\newcommand{\Hbf}{{ H}}
\newcommand{\I}{{ I}}
\newcommand{\R}{{ R}}
\newcommand{\Sbf}{{S}}
\newcommand{\Sigmabf}{{\Sigma}}

\newremark{remark}{Remark}
\newtheorem{theorem}{Theorem}
\newtheorem{corollary}{Corollary}
\newtheorem{lemma}{Lemma}
\newcommand{\cal}{\mathcal}
\newcommand{\fraca}[2]{{#1}/{#2}}
\makeatother

\begin{document}
\begin{frontmatter}

\title{Empirical likelihood for single-index varying-coefficient models}
\runtitle{Single-index varying-coefficient models}

\begin{aug}
\author[1]{\fnms{Liugen} \snm{Xue}\corref{}\thanksref{1}\ead[label=e1]{lgxue@bjut.edu.cn}} \and
\author[2,3]{\fnms{Qihua} \snm{Wang}\thanksref{2,3}\ead[label=e2]{qhwang@amss.ac.cn}}
\address[1]{College of Applied Sciences, Beijing University of
Technology, Beijing 100124, China.\\\printead{e1}}
\address[2]{Academy of Mathematics and Systems Science, Chinese
Academy of Science, Beijing 100080, China}
\address[3]{School of Mathematics and Statistics, Yunnan University,
Kunming 650091, China.\\\printead{e2}}
\end{aug}

\received{\smonth{8} \syear{2010}}

%
\begin{abstract}
In this paper, we develop statistical
inference techniques for the unknown coefficient functions and
single-index parameters in single-index varying-coefficient
models. We first estimate the nonparametric component via the
local linear fitting, then construct an estimated empirical
likelihood ratio function and hence obtain a maximum empirical
likelihood estimator for the parametric component. Our estimator
for parametric component is asymptotically  efficient,  and the estimator of nonparametric component
has an optimal convergence rate. Our results provide ways to
construct the confidence region for the involved unknown
parameter. We also develop an adjusted empirical likelihood ratio
for constructing the confidence regions of parameters of interest.
A simulation study is conducted to evaluate the finite sample
behaviors of the proposed methods.
\end{abstract}

%

\end{frontmatter}\vspace*{-2pt}

\section{Introduction}\label{sec1}

Consider a single-index varying-coefficient model of the form
%
\begin{equation}
Y={\rm g}_0^T(\beta_0^TX)Z+\varepsilon,
\label{eq1.1}
\end{equation}
where $(X, Z)\in R^p\times R^q$ is a vector of covariates, $Y$ is
the response variable, $\beta_0$ is an $p\times1$ vector of unknown
parameters, ${\rm g}_0(\cdot)$ is an $q\times1$ vector of unknown
functions and~$\varepsilon$ is a random error with mean 0 and finite
variance $\sigma^2$. Assume that $\varepsilon$ and $(X,Z)$ are
independent. For the sake of identifiability, it is often assumed
that $\|\beta_0 \|=1$, and the first non-zero element is positive,
where $\|\cdot\|$ denotes the Euclidean metric.

Model (\ref{eq1.1}) includes a class of important statistical models. For
example, if $q=1$ and $Z=1$, (\ref{eq1.1}) reduces to the single-index model
(see, e.g., H\"ardle, Hall and Ichimura~\cite{HarHalIch93},
Weisberg
and Welsh \cite{WeiWel94}, Zhu and Fang \cite{ZhuFan96},
Chiou and M\"uller \cite{ChiMul98}, Hristache, Juditsky and Spokoiny
\cite{HriJudSpo01}, Xue and Zhu \cite{XueZhu06}). If $p=1$ and
$\beta_0=1$, (\ref{eq1.1}) is the varying-coefficient model (see, e.g.,
Chen and Tsay \cite{CheTsa93}, Hastie and Tibshirani \cite{HasTib93},
Wu, Chiang and Hoover~\cite{WuChiHoo98}, Fan and Zhang~\cite
{FanZha99}, Cai, Fan and Li \cite{CaiFanLi00},
Cai, Fan and Yao \cite{CaiFanYao00}, Xue and Zhu \cite{XueZhu07N1}).
If the last component of $\beta_0$ to be non-zero and $Z=(1,X^*T)^T$ where $X^*$ is the remaining vector of
$X$ with its $p$th component deleted, (\ref{eq1.1}) becomes the
adaptive varying-coefficient linear model (see, e.g., Fan, Yao and Cai
\cite{FanYaoCai03}, Lu, Tj\o stheim and Yao~\cite{LuTjsYao07}).

Model (\ref{eq1.1}) is easily interpreted in real applications because it
has the features of the single-index model and the
varying-coefficient model. In addition, model (\ref{eq1.1}) may include
cross-product terms of some components of $X$ and $Z$. Hence it
has considerable flexibility to cater for complex multivariate
nonlinear structure. Xia and Li \cite{XiaLi99} investigated a class of
single-index coefficient regression models, which include model
(\ref{eq1.1}) as a~special example. When it is used as a~nonparametric
time series model, Xia and Li \cite{XiaLi99} obtained the estimator of
${{\rm g}}(\cdot)$ by kernel smoothing and then derived the estimator
of $\beta_0$ by the least squares method and proved that the
corresponding estimators are consistent and asymptotically normal.

In this paper, we develop statistical inference techniques of
${{\rm g}}_0(\cdot)$ and $\beta_0$ with independent observations of
$(Y, X, Z)$. We can construct an empirical likelihood ratio
function for $\beta_0$ by assuming ${\rm g}_0(\cdot)$ and its
derivative to be known functions. In practice, however, they are
unknown, and hence the empirical likelihood ratio function cannot
be used to make inference on $\beta$. This motivates us to
estimate the unknown ${\rm g}_0(\cdot)$ and $\dot{{\rm g}}_0(\cdot)$
via the
local linear smoother, and then obtain an estimated empirical
likelihood ratio of $\beta_0$. The estimated empirical
log-likelihood ratio is asymptotically distributed as a weighted
sum of independent $\chi_1^2$ variables with unknown weights. This
result cannot be applied directly to construct confidence region
for $\beta_0$. To solve this issue, two methods may be used (see
Wang and Rao \cite{WanRao02}). The first method is to estimate the unknown
weights consistently so that the distribution of the estimated
weighted sum of chi-squared variables can be estimated from the
data. The second method is to adjust the estimated empirical
log-likelihood ratio so that the resulting adjusted empirical
log-likelihood ratio is asymptotically chi-squared. Also, we obtain
a maximum empirical likelihood estimator of $\beta_0$, by
maximizing the estimated empirical likelihood ratio function, and
investigate its asymptotic property. In addition, we obtain the
convergence rate of the estimator of~$\sigma^2$ and define the
consistent estimator of asymptotic variance; this allows us to
construct a~confidence region for $\beta_0$.

Comparing with the existing methods, our estimating method has the
following advantage: The asymptotic variance of our estimator for
$\beta_0$ is the same as  those of H\"ardle \textit{et al.} \cite
{HarHalIch93} and Xia and Li \cite{XiaLi99} when the model reduces to
the single-index model;
this shows that our estimator for $\beta_0$ is the same efficient as than
those of H\"ardle \textit{et al.} \cite{HarHalIch93} and Xia and Li \cite
{XiaLi99}.
The difference between the proposed estimating approaches and the
existing estimating approaches is that we use an empirical likelihood
ratio to define the estimator of $\beta_0$ while the existing work
uses the least squares techniques (see, e.g., H\"ardle \textit{et al.}
\cite{HarHalIch93}, Xia and Li \cite{XiaLi99}). Also, we develop an
empirical likelihood
inference for constructing a confidence region of $\beta$. The
empirical likelihood method, introduced by Owen\vadjust{\goodbreak} \cite{Owe88}, has many
advantages for constructing confidence regions or intervals. For
example, it does not impose prior constraints on region shape, and
it does not require the construction of a pivotal quantity. The
empirical likelihood has been studied by many authors. The related
works are Wang and Rao \cite{WanRao02}, Wang, Linton and H\"{a}rdle
\cite{WanLinHar04},
Xue and Zhu \cite{XueZhu06,XueZhu07N1,XueZhu07N2},
Zhu and Xue \cite{ZhuXue06}, Qin and Zhang \cite{QinZha07}, Stute,
Xue and Zhu \cite{StuXueZhu07},
Xue \cite{Xue09N2,Xue09N1}, Wang and Xue \cite{WanXue11}, among others.

The rest of the paper is organized as follows. In Section~\ref{sec2}, we
define an estimated empirical likelihood ratio, and then obtain a
maximum empirical likelihood estimator of~$\beta_0$ by maximizing
the empirical likelihood ratio function; the asymptotic properties
of the proposed estimators are also investigated. In Section~\ref{sec3}, we
define an adjusted empirical log-likelihood and derive its
asymptotic distribution. Section~\ref{sec4} reports a simulation study.
Proofs of theorems are relegated
to the \hyperref[appm]{Appendix}. It should be pointed that some special techniques
are used in the proofs.

\section{Estimated empirical likelihood}\label{sec2}

\subsection{Methodology}\label{sec2.1}

Suppose that $ \{(Y_i,X_i, Z_i);1\leq i\leq n \}$ is an
independent and identically distributed (i.i.d.) sample from (\ref{eq1.1}),
that is
\[
Y_i={\rm g}_0^T(\beta_0^TX_i)Z_i+\varepsilon_i,   \qquad  i=1,\ldots,n,
\]
where $\varepsilon_i$s are i.i.d. random errors with mean 0 and
finite variance $\sigma^2$. Assume that
$ \{\varepsilon_i;1\leq i\leq n \}$ are independent of
$ \{(X_i,Z_i);1\leq i\leq n \}$.

To construct an empirical likelihood ratio function for $\beta_0$,
we introduce an auxiliary random vector
%
\begin{equation}
\eta_i(\beta)=\{Y_i-{\rm g}_0^T(\beta^TX_i)Z_i\}\dot{{\rm
g}}_0^T(\beta^TX_i)Z_iX_iw(\beta^TX_i),
\label{eq2.1}
\end{equation}
where $\dot{{\rm g}}_0(\cdot)$ stands for the derivative of the
function vector ${\rm g}_0(\cdot)$, and $w(\cdot)$ is a bounded weight
function with a bounded support ${\cal U}_w$, which is introduced
to control the boundary effect in the estimations of
${{\rm g}}_0(\cdot)$ and $\dot{{\rm g}}_0(\cdot)$. To convenience, we take that $w(\cdot)$ is the indicator function of the set
${\cal U}_w$. Note that
$E\{\eta_i(\beta)\}=0$ if $\beta=\beta_0$. Hence, the problem of
testing whether $\beta$ is the true parameter is equivalent to
testing whether $E\{\eta_i(\beta)\}=0$ for $i=1,2,\ldots,n$.
By Owen \cite{Owe88}, this can be done by
using the empirical likelihood. That is, we can define the profile
empirical likelihood ratio function
\[
{L}_n(\beta)=\max \Biggl\{\prod_{i=1}^n(np_i) \bigg| p_i\geq0,
\sum_{i=1}^n p_i=1, \sum_{i=1}^np_i{\eta}_i(\beta)=0 \Biggr\}.
\]
It can be shown that $-2\log L_n(\beta_0)$ is asymptotically
chi-squared with $p$ degrees of freedom. However, $L_n(\beta)$
cannot be directly used to make statistical inference on $\beta_0$
because $L_n(\beta)$ contains the unknowns ${\rm g}_0(\cdot)$ and
$\dot{{\rm g}}_0(\cdot)$. A natural way is to replace ${\rm
g}_0(\cdot)$ and
$\dot{{\rm g}}_0(\cdot)$ in ${L}_n(\beta)$ by their estimators and define
an estimated empirical likelihood function. In this paper, we
estimate the vector functions ${\rm g}_0(\cdot)$ and $\dot{{\rm
g}}_0(\cdot)$
via the local linear regression technique (see, e.g., Fan and Gijbels
\cite{FanGij96}). The local linear estimators for ${\rm g}_0(u)$ and
$\dot{{\rm g}}_0(u)$ are defined as $\hat{{\rm g}}(u;\beta
_0)=\hat{\rm a}$
and $\hat{\dot{{\rm g}}}(u;\beta_0)=\hat{\rm b}$ at the fixed point
$\beta_0$, where $\hat{\rm a}$ and $\hat{\rm b}$ minimize the sum of
weighted squares
\[
\sum_{i=1}^n [Y_i-\{{\rm a} + {\rm b}
(\beta_0^TX_i-u)\}^TZ_{i} ]^2K_h(\beta_0^TX_i-u),
\]
where $K_h(\cdot)=h^{-1}K(\cdot/h)$, $K(\cdot)$ is a kernel
function, and $h=h_n$ is a bandwidth sequence that decreases to 0
as $n$ increases to $\infty$. It follows from the least squares
theory that
\[
 (\hat{{\rm g}}^T (u;\beta_0),h\hat{\dot{{\rm g}}}{}^T
(u;\beta_0) )^T
={\Sbf}_n^{-1}(u;\beta_0)\xi_n(u;\beta_0),
\]
where
%
\[
{\Sbf}_n(u;\beta_0)=\left(
\begin{array}{c@{ \quad }c}
{\Sbf}_{n,0}(u;\beta_0) & {\Sbf}_{n,1}(u;\beta_0) \\
{\Sbf}_{n,1}(u;\beta_0) & {\Sbf}_{n,2}(u;\beta_0)
\end{array}
\right)
 \quad \mbox{and} \quad
\xi_n(u;\beta_0)=\left(
\begin{array}{c}
\xi_{n,0}(u;\beta_0) \\
\xi_{n,1}(u;\beta_0)
\end{array}
\right)
\]
with\vspace*{-1pt}
\[
{\Sbf}_{n,j}(u;\beta_0)=\frac{1}{n}\sum_{i=1}^nZ_iZ_i^T \biggl(\frac
{\beta_0^TX_i-u}{h} \biggr)^jK_h(\beta_0^TX_i-u)
\]
and\vspace*{-1pt}
\[
\xi_{n,j}(u;\beta_0)=\frac{1}{n}\sum_{i=1}^nZ_iY_i \biggl(\frac
{\beta_0^TX_i-u}{h} \biggr)^jK_h(\beta_0^TX_i-u).
\]
%

Since the convergence rate of the estimator of $\dot{g}_0'(u)$ is slower than that of the estimator of ${g}_0(u)$ if the same bandwidth is used,
this leads to a slower convergence rate for the estimator $\hat\beta$ of $\beta_0$ than
 $\sqrt n$. To increase the convergence rate of the estimator of $\dot{g}_0'(u)$, we introduce the another bandwidth $h_1$ to replace $h$ in $\hat{\dot{\rm g}}(u;\beta)$, and define as $\hat{\dot{\rm g}}_{h_1}(u;\beta)$.

Let $\hat{\eta}_i(\beta)$ be $\eta_i(\beta)$, with ${\rm
g}_0(\beta^TX_i)$
and $\dot{{\rm g}}_0(\beta^TX_i)$ replaced by
$\hat{{\rm g}}(\beta^TX_i;\beta)$ and $\hat{\dot{{\rm g}}}_{h_1}(\beta
^TX_i;\beta)$,
respectively, for $i=1,\ldots,n$. Then an estimated empirical
likelihood ratio function is defined by
\[
\hat{L}(\beta)=\max \Biggl\{\prod_{i=1}^n(np_i) \bigg| p_i\geq0,
\sum_{i=1}^n p_i=1, \sum_{i=1}^np_i\hat{\eta}_i(\beta)=0 \Biggr\}.
\]
By the Lagrange multiplier method, $\log\hat{L}(\beta)$ can be
represented as
%
\begin{equation}
\log\hat{L}(\beta)=-\sum_{i=1}^n\log\bigl (1+\lambda^T\hat{\eta
}_i(\beta) \bigr),
\label{eq2.2}
\end{equation}
where $\lambda$ is determined by
%
\begin{equation}
\frac{1}{n}\sum_{i=1}^n\frac{\hat{\eta}_i(\beta)}{1+\lambda
^T\hat{\eta}_i(\beta)}=0.
\label{eq2.3}
\end{equation}

Let ${\cal B}=\{\beta\in R^p\dvt \|\beta\|=1$, and the first
non-zero element is positive. Then $\beta_0$ is an inner point
of the set ${\cal B}$. Therefore we need only search for $\beta_0$
over ${\cal B}$. A maximum empirical likelihood estimator for
$\beta_0$ is given by
%
\begin{equation}
\hat{\beta}=\arg\sup_{\beta\in{\cal B}}\hat{L}(\beta).
\label{eq2.4}
\end{equation}
With $\hat{\beta}$, we define the estimate of ${\rm g}(u)$ by
$
\hat{{\rm g}}(u)=\hat{{\rm g}}(u, \hat{\beta}),
$
and the estimate of $\sigma^2$ by
%
\begin{equation}
\hat{\sigma}^2=\frac{1}{n}\sum_{i=1}^n\{Y_i-\hat{{\rm g}}{}^T(\hat
{\beta}^TX_i;\hat{\beta})Z_i\}^2.
\label{eq2.5}
\end{equation}

It is well known that if $\beta$ is known, the optimal bandwidth $h$
for $\hat{{\rm g}}(u)$ is of order $\mathrm{O}(n^{-1/5})$. However, if $\beta$ is
unknown, in order to ensure that the estimator $\hat{\beta}$ is
root-$n$ consistent, the bandwidth $h$ should be smaller than
$\mathrm{O}(n^{-1/5})$, if we only assume ${{\rm g}}(\cdot)$ are second-order
differentiable (see Theorem~\ref{theo2} below). Note that once the estimator
$\hat{\beta}$ is available, an optimal bandwidth of order
$\mathrm{O}(n^{-1/5})$ can be used in the final estimator for ${{\rm g}}(\cdot)$.

%

\subsection{Asymptotic properties}\label{sec2.2}

In order to obtain the asymptotic behaviors of our estimators, we
first give the following conditions:
\begin{enumerate}[(C6)]
\item[(C1)]The density function of $\beta^TX$, $f(u)$, is bounded
away from zero for $u\in{\cal U}_w$ and $\beta$ near $\beta_0$,
and satisfies the Lipschitz condition of order 1 on ${\cal U}_w$,
where ${\cal U}_w$ is the support of $w(u)$.

\item[(C2)]The functions $g_{j}(u)$, $1\leq j\leq q$, have
continuous second derivatives on ${\cal U}_w$, where $g_{j}(u)$
are the $j$th components of ${\rm g}_0(u)$.

\item[(C3)]$E(\|X\|^6)<\infty$, $E(\|Z\|^6)<\infty$ and
$E(|\varepsilon|^6)<\infty$.

\item[(C4)]$nh^2/\log^2{n}\rightarrow\infty$, $nh^4\log{n}\rightarrow0$; $nhh_1^3/\log^2{n}\rightarrow\infty$, $nh_1^5=\mathrm{O}(1)$.

\item[(C5)]The kernel $K(\cdot)$ is a symmetric probability
density function with a bounded support and satisfies the
Lipschitz condition of order 1 and $\int u^2K(u)\,\mathrm{d}u\neq0$.

\item[(C6)]The matrix ${\D}(u)=E(ZZ^T|\beta_0^TX=u)$ is positive
definite, and each entry of ${\D}(u)$ and
${\C}(u)=E(VZ^T|\beta_0^TX=u)$ satisfies the Lipschitz condition of
order 1 on ${\cal U}_w$, where
$V=X\dot{{\rm g}}_0^T(\beta_0^TX)Zw(\beta_0^TX)$, and ${\cal U}_w$ is
defined in (C1).

\item[(C7)]The matrices ${B}(\beta_0)=E(VV^T)$ and ${B}_*(\beta_0)=
{B}(\beta_0)-E\{C(\beta_0^TX)\dot{\rm g}_0(\beta_0^TX)E(X^T|\allowbreak \beta_0^TX)\}$ are positive definite, where $V$ is defined in (C6).
\end{enumerate}

\begin{remark}     Condition {\rm(C1)} is used to bound the
density function of $\beta^TX$ away from zero. This ensures that
the denominators of $\hat{{\rm g}}(u;\beta)$ and
$\hat{\dot{{\rm g}}}(u;\beta)$ are, in probability one, bounded away
from 0 for $u\in{\cal U}_w$. The second derivatives in {\rm(C2)}
are standard smoothness conditions. {\rm(C3)--(C5)} are necessary
conditions for the asymptotic normality or the uniform consistency
of the estimators. Conditions {\rm(C6)} and {\rm(C7)} ensure that
the asymptotic variance for the estimator of $\beta_0$ exists.

Let ${\cal B}_n=\{\beta\in{\cal B}\dvt \|\beta-\beta_0\|\leq
c_0n^{-1/2}\}$ for some positive constant $c_0$. This is motivated
by the fact that, since we anticipate that $\hat{\beta}$ is
root-$n$ consistent, we should look for a~maximum of
$\hat{L}(\beta)$ which involves $\beta$ distant from $\beta_0$ by
order $n^{-1/2}$. Similar restrictions were also made by H\"ardle,
Hall and Ichimura \cite{HarHalIch93}, Xia and Li \cite{XiaLi99} and
Wang and Xue~\cite{WanXue11}.

The following theorem shows that $-2\log\hat{L}(\beta_0)$ is
asymptotically distributed as a~weighted sum of independent $\chi_1^2$
variables.
\end{remark}

\begin{theorem}\label{theo1}  Suppose that conditions {\rm(C1)--(C7)}
hold. Then
\[
-2\log\hat{L}(\beta_0)\stackrel{D}{\longrightarrow}w_1\chi
_{1,1}^2+\cdots +w_p\chi_{1,p}^2,
\]
where $\stackrel{D}{\longrightarrow}$ represents convergence in
distribution, $\chi_{1,1}^2, \ldots,\chi_{1,p}^2$ are independent
$\chi_1^2$ variables and the weights $w_j$, for $1\leq j\leq p$,
are the eigenvalues of
${\G}(\beta_0)={\B}^{-1}(\beta_0){\A}(\beta_0)$. Here
${\B}(\beta_0)$ is defined in condition {\rm(C7)},
%
\begin{equation}
{\A}(\beta_0)={\B}(\beta_0)-E\{{\C}(\beta_0^TX){\D}^{-1}(\beta
_0^TX){\C}{}^T(\beta_0^TX)\},
\label{eq2.6}
\end{equation}
and  ${\C}(u)$ and ${\D}(u)$ are defined in condition {\rm(C6)}.
\end{theorem}

To apply Theorem~\ref{theo1} to construct a confidence region or interval
for $\beta_0$, we need to consistently estimate the unknown
weights $w_j$. By the ``plug-in'' method, ${\A}(\beta_0)$ and~${\B}(\beta_0)$ can be consistently estimated by
%
\begin{equation}
\hat{{\A}}(\hat{\beta})=\frac{1}{n}\sum_{i=1}^n \{\hat
{V}_i\hat{V}_i^T
-
\hat{{\C}}(\hat{\beta}^TX_i)\hat{{\D}}^{-1}(\hat{\beta
}^TX_i)\hat{{\C}}{}^T(\hat{\beta}^TX_i) \}
\label{eq2.7}
\end{equation}
and
%
\begin{equation}
\hat{{\B}}(\hat{\beta})=\frac{1}{n}\sum_{i=1}^n\hat{V}_i\hat{V}_i^T,
\label{eq2.8}
\end{equation}
respectively, where $\hat{\beta}$ is the maximum empirical
likelihood estimator of $\beta_0$ defined by\vspace*{1pt}~(\ref{eq2.4}),
$\hat{V}_i=X_i\hat{\dot{{\rm g}}}{}^T(\hat{\beta}^TX_i;\hat{\beta
})Z_iw(\hat{\beta}^TX_i)$,
$\hat{{\C}}(\cdot)=\sum_{i=1}^nW_{ni}(\cdot)\hat{V}_iZ_i^T$ and
$\hat{{\D}}(\cdot)=\break\sum_{i=1}^nW_{ni}(\cdot)Z_iZ_i^T$ with
\[
W_{ni}(\cdot)=K_1 \biggl(\frac{\hat{\beta}^TX_i-\cdot}{b_n}
\biggr) \bigg/\sum_{k=1}^nK_1 \biggl(\frac{\hat{\beta}^TX_k-\cdot
}{b_n} \biggr),
\]
where $K_1(\cdot)$ is a kernel function, and $b_n$ is a bandwidth
with $0<b_n\rightarrow0$.

This implies that the eigenvalues of
$\hat{\G}(\hat{\beta})=\hat{\B}^{-1}(\hat{\beta})\hat{\A
}(\hat{\beta})$,
say $\hat{w}_j$, consistently estimate~$w_j$ for $j=1,\ldots,p$.
Let $\hat{c}_{1-\alpha}$ be the $1-\alpha$ quantile of the
conditional distribution of the weighted sum $\hat{s} =
\hat{w}_1\chi_{1,1}^2+\cdots +\hat{w}_p\chi_{1,p}^2$ given the data.
Then an approximate $1-\alpha$ confidence region for $\beta_0$ can
be defined as
\[
{\cal R}_{\rm eel}(\alpha) =  \{\beta\in{\cal B}\dvt
-2\log\hat{L}(\beta)\leq\hat{c}_{1-\alpha} \}.
\]
In practice, the conditional distribution of the weighted sum
$\hat{s}$, given the sample $\{(Y_i,X_i, Z_i),1\leq i\leq n\}$, can be
calculated using Monte Carlo simulations by repeatedly generating
independent samples $\chi_{1,1}^2,\ldots,\chi_{1,p}^2$ from the
$\chi_1^2$ distribution.

The following theorem states an interesting result about $\hat
\beta$. The asymptotic variance of $\hat\beta$ is smaller than that
of H\"ardle \textit{et al.} \cite{HarHalIch93} when our model reduces to
a single-index
model.\looseness=-1

\begin{theorem}\label{theo2}  Suppose that conditions {\rm(C1)--(C7)}
hold. Then
\[
\sqrt{n} (\hat{\beta}-\beta_0 )\stackrel
{D}{\longrightarrow}N (0,\sigma^{2}{\B_*^{-1}}(\beta_0){\A
}(\beta_0){\B_*^{-1}}(\beta_0) ),
\]
  where ${\B_*}(\beta_0)$ and ${\A}(\beta_0)$ are defined in
condition {\rm(C7)} and {\rm(\ref{eq2.6})}, respectively.
\end{theorem}

In model (\ref{eq1.1}), if $q=1$ and $Z=1$, then (\ref{eq1.1}) reduces to the
single-index model. By Theorem~\ref{theo2}, we derive the following result.

\begin{corollary}\label{cor1}  Suppose that the conditions of Theorem~\ref{theo2}
hold. If $q=1$ and $Z=1$ in model (\ref{eq1.1}), then
\[
\sqrt{n}(\hat{\beta}-\beta_0)\stackrel{D}{\longrightarrow}N
(0,\sigma^{2}{\A}_1^-(\beta_0)),
\]
  where
${\A}_1(\beta_0)=E[\{X-E(X|\beta_0^TX)\}\{X-E(X|\beta_0^TX)\}^T\dot
{{\rm g}}_0^2(\beta_0^TX)w(\beta_0^TX)]$
and $A_1^{-}$ represents a generalized inverse of the matrix $A_1^{-}$.
\end{corollary}

Corollary~\ref{cor1} is the same as the results of H\"{a}rdle et al. \cite{HarHalIch93} and Xia and Li \cite{XiaLi99}
 for the single-index model.

For the estimator of the variance of error, $\hat{\sigma}^2$, we
have the following result.

\begin{theorem}\label{theo3}  Suppose that conditions   \textup{(C1)--(C7)}
hold. Then,
\[
\hat{\sigma}^2-\sigma^2 = \mathrm{O}_P (n^{-1/2} ).
\]
\end{theorem}

To apply Theorem~\ref{theo2} to construction of the confidence region of $\beta_0$, we use the estimators
$\hat{\sigma}^2$ and $\hat{A}(\hat{\beta})$ defined in (\ref{eq2.5}) and (\ref{eq2.7}), and define the estimator of $B_*(\beta_0)$ as follows
\[
\hat{B}_*(\hat{\beta})=\frac{1}{n}\sum_{i=1}^n\{\hat{V}_i\hat{V}_i^T-\hat{C}(\hat{\beta}^TX_i)\hat{\dot{g}}(\hat{\beta}^TX_i;\hat\beta)\hat{\mu}^T(\hat{\beta}^TX_i) \},
\]
where $\hat{\mu}(\cdot)=\sum_{i=1}^nW_{ni}(\cdot)X_i$ is the estimator of $\mu(u)=E(X|\beta_0^TX=u)$.
It can be shown that $\hat{A}(\hat{\beta})\stackrel{P}{\longrightarrow}{A}(\beta_0)$ and $\hat{B}_*(\hat{\beta})\stackrel{P}{\longrightarrow}{B}_*(\beta_0)$,
where $\stackrel{P}{\longrightarrow}$ denotes convergence in
probability. By Theorems~\ref{theo3} and~\ref{theo4}, we have
\[
 \{\hat{\sigma}^2{\hat{\B}_*^{-1}}(\hat{\beta})\hat{{\A
}}(\hat{\beta}){\hat{\B}_*^{-1}}(\hat{\beta}) \}^{-1/2}
\sqrt{n}(\hat{\beta}-\beta_0)
\stackrel{D}{\longrightarrow}N(0,{\I}_p).
\]
Using Theorem 10.2d in Arnold \cite{Arn81}, we obtain
\[
(\hat{\beta}-\beta_0)^T
 \{n^{-1}\hat{\sigma}^2{\hat{\B}_*^{-1}}(\hat{\beta})\hat{{\A
}}(\hat{\beta}){\hat{\B}_*^{-1}}(\hat{\beta}) \}^{-}
(\hat{\beta}-\beta_0)
\stackrel{D}{\longrightarrow}\chi_p^2.\vadjust{\goodbreak}
\]
%

Let $\chi_p^2(1-\alpha)$ be the $1-\alpha$ quantile of $\chi_p^2$
for $0<\alpha<1$. Then
\[
 \{\beta\dvt (\hat{\beta}-\beta)^T
 (n^{-1}\hat{\sigma}^2{\hat{\B}_*^{-1}}(\hat{\beta})\hat{{\A
}}(\hat{\beta})
{\hat{\B}_*^{-1}}(\hat{\beta}) )^{-}
(\hat{\beta}-\beta)\le\chi_p^2(1-\alpha) \}
\]
gives an approximate $1-\alpha$ confidence region for $\beta_0$.

\section{Adjusted empirical likelihood}\label{sec3}

In addition to the above, direct way of approximating the
asymptotic distributions, we can also consider the following
alternative. The alternative is motivated by the results of Rao and
Scott \cite{RaoSco81}. By Rao and Scott \cite{RaoSco81} the
distribution of
$\rho(\beta_0)\sum_{i=1}^p w_i\chi_{1,i}^2$ can be approximated by
$\chi_p^2$, where $\rho(\beta_0)=p/\operatorname{tr}\{{\G}(\beta_0)\}$. Let
$\hat\rho(\hat\beta)=p/\operatorname{tr}\{\hat{\G}(\hat\beta)\}$ with
$\hat{\G}(\hat{\beta})=\hat{\A}^{1/2}(\hat{\beta})\hat{\B
}^{-1}(\hat{\beta})\hat{\A}^{1/2}(\hat{\beta})$,
where $\hat{\A}(\hat{\beta})$ and $\hat{\B}(\hat{\beta})$ are
defined in (\ref{eq2.7}) and (\ref{eq2.8}). Invoking Theorem~\ref{theo1} and the consistency
of $\hat{\G}(\hat{\beta})$, the asymptotic distribution of
$\hat{\rho}(\hat{\beta})\{-2\log\hat{L}(\beta)\}$ can be
approximated by $\chi_p^2$. Clearly, $\hat\beta$ in
$\hat\rho(\cdot)$ can be replaced by $\beta$. Therefore, an
improved Rao--Scott adjusted empirical log-likelihood can be defined
as
\[
\tilde{l}(\beta)=\hat{\rho}(\beta)\{-2\log\hat{L}(\beta)\}.
\]
However, the accuracy of this approximation still depends on the
values of the $w_i$s. Now, we propose another adjusted empirical
log-likelihood, whose asymptotic distribution is chi-squared with $p$
degrees of freedom. The adjustment technique is developed by Wang and
Rao \cite{WanRao02} by using an approximate result in Rao and Scott
\cite{RaoSco81}. Note that $\hat{\rho}(\beta)$ can be written as
\[
\hat{\rho}(\beta)= \frac{\operatorname{tr}\{\hat{\A}^{-}(\beta)\hat{\A}(\beta)\}}
{\operatorname{tr}\{\hat{\B}^{-1}(\beta)\hat{\A}(\beta)\}}.
\]
By examining the asymptotic expansion of $-2\log\hat{L}(\beta)$,
which is specified in the proof of Theorem~\ref{theo4} below, we define an
adjustment factor
\[
\hat{r}({\beta})= \frac{\operatorname{tr}\{\hat{\A}^{-}(\beta)\hat{\Sigmabf}(\beta)\}}
{ \operatorname{tr}\{\hat{\B}^{-1}(\beta)\hat{\Sigmabf}(\beta)\}},
\]
by replacing ${\hat A}(\beta)$ in ${\hat\rho}(\beta)$ by ${\hat
{\Sigmabf}}(\beta) $, where $\hat{\Sigmabf}(\beta)
= \{\sum_{i=1}^n\hat{\eta}_i(\beta) \} \{\sum
_{i=1}^n\hat{\eta}_i(\beta) \}^T$.
The adjusted empirical log-likelihood ratio is defined by
%
\begin{equation}
\hat{l}_{\rm ael}(\beta)=\hat{r}(\beta)\{-2\log\hat{L}(\beta)\},
\label{eq3.1}
\end{equation}
where $\log\hat{L}(\beta)$ is defined in (\ref{eq2.2}).

\begin{theorem}\label{theo4}  Suppose that conditions {\rm(C1)--(C6)} hold.
Then
$
\hat{l}_{\rm ael}(\beta_0)\stackrel{D}{\longrightarrow}\chi_p^2.
$
\end{theorem}

According to Theorem~\ref{theo4}, $\hat{l}_{\rm ael}(\beta)$ can be used to
construct an approximate confidence region for $\beta_0$. Let
\[
{\cal R}_{\rm ael}(\alpha) =  \{\beta\in{\cal B}\dvt \hat{l}_{\rm
ael}(\beta)\leq\chi_p^2(1-\alpha) \}.
\]
Then, ${\cal R}_{\rm ael}(\alpha)$ gives a confidence region for
$\beta_0$ with asymptotically correct coverage probability
$1-\alpha$.

\section{Numerical results}\label{sec4}


\subsection{Bandwidth selection}\label{sec4.1}

Various existing bandwidth selection techniques
for nonparametric regression, such as the cross-validation and
generalized cross-validation, can be adapted for the estimation~$\hat{{\rm g}}(\cdot)$. But we, in our simulation, use the modified
multi-fold cross-validation (MMCV) criterion proposed by Cai, Fan and
Yao \cite{CaiFanYao00} to select the optimal bandwidth because the
algorithm is simple and quick. Let $m$ and $Q$ be two given
positive integers and $n>mQ$. The basic idea is first to use $Q$
sub-series of lengths $n-km$ $(k=1,\ldots,Q)$ to estimate the
unknown coefficient functions and then to compute the one-step
forecasting error of the next section of the sample of lengths $m$
based on the estimated models. More precisely, we choose $h$ which
minimizes
%
\begin{equation}
\operatorname{AMS}(h) = \sum_{k=1}^Q\operatorname{AMS}_k(h),
\label{eq4.1}
\end{equation}
where, for $k=1,\ldots,Q$,
\[
\operatorname{AMS}_k(h)=\frac{1}{m}\sum_{i=n-km+1}^{n-km+m}\Biggl \{Y_i
-\sum_{j=1}^q\hat{g}_{j,k}(U_i)Z_{ij} \Biggr\}^2,
\]
and $\{\hat{g}_{j,k}(\cdot)\}$ are computed from the sample
$\{(Y_i,U_i,Z_i),1\leq i\leq n-km\}$ with bandwidth equal
$h(\frac{n}{n-km})^{1/5}$. Note that for different sample size, we
re-scale bandwidth according to its optimal rate, that is, $h\propto
n^{-1/5}$. Since the selected bandwidth does not depend critically
on the choice of $m$ and $Q$, to computation expediency, we take
$m=[0.1n]$ and $Q=4$ in our simulation.

Let $h_{\rm opt}$ be the bandwidth obtained by minimizing (\ref{eq4.1})
with respect to $h > 0$; that is, $h_{\rm
opt}=\inf_{h>0}\operatorname{AMS}(h)$. Then $h_{\rm opt}$ is the optimal
bandwidth for estimating $\hat{{\rm g}}(\cdot)$. When calculating the
empirical likelihood ratios and estimator of $\beta_0$, we use the
approximation bandwidth
\[
h=h_{\rm opt}n^{-1/20}(\log n)^{-1/2}, \qquad h_1=h_{\rm opt},
\]
because this insures that the required bandwidth has correct order
of magnitude for the optimal asymptotic performance (see, e.g.,
Carroll \textit{et al.} \cite{Caretal97}), and the bandwidth $\hat{h}$ satisfies
condition~(C4).

\subsection{Simulation study}\label{sec4.2}

We now examine the performance of the procedures described in
Sections~\ref{sec2} and~\ref{sec3}. Consider the regression model
%
\begin{equation}
Y_i = g_0(\beta_0^TX_i) + g_1(\beta_0^TX_i)Z_{i1} +
g_2(\beta_0^TX_i)Z_{i2} + \varepsilon_i,
\label{eq4.2}
\end{equation}
where $\beta_0=(1/\sqrt{5},2/\sqrt{5})^T$ and the $\varepsilon_i$s
are independent $N(0,0.8^2)$ random variables. The sample
$\{X_i=(X_{i1},X_{i2})^T;1\leq i\leq n\}$ was generated from a
bivariate uniform distribution on $[-1,1]^2$ with independent
components, $\{Z_i=(Z_{i1},Z_{i2})^T;1\leq i\leq n\}$ was
generated from a bivariate normal distribution $N(0,\Sigma)$ with
$\operatorname{var}(Z_{i1})=\operatorname{var}(Z_{i2})=1$ and the correlation
coefficient between $Z_{i1}$ and $Z_{i2}$ is $\rho=0.6$. In
model (\ref{eq4.2}), the coefficient functions are $g_0(u)=12\exp(-2u^2)$,
$g_1(u)=10u^2$ and $g_2(u)=16\sin(\uppi u)$.

For the smoother, we used a local linear smoother with a
Epanechnikov kernel $K(u)=0.75(1-u^2)_+$ with a MMCV bandwidth
throughout all smoothing steps. We take the weight function
$w(u)=I_{[-3/\sqrt{5},3/\sqrt{5}]}(u)$. The sample size for the
simulated data is $100$,\vspace*{2pt} and the run is 500 times in all
simulations.

\begin{figure}

\includegraphics{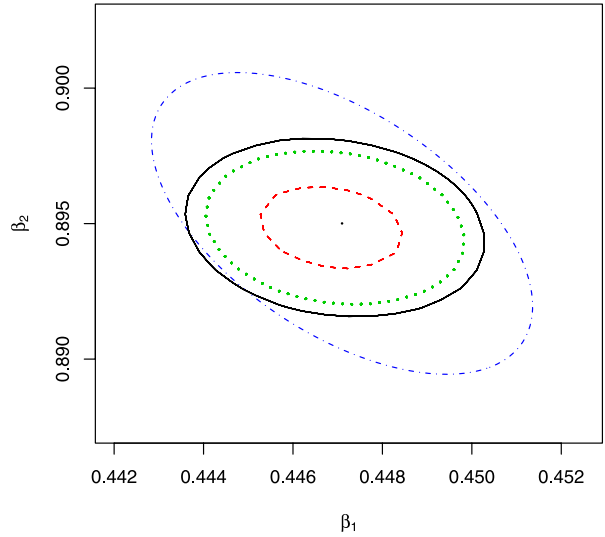}

\caption{Averages of 95\%
confidence regions of $(\beta_1,\beta_2)$, based on EEL (solid
curve), AEL (dashed curve), IRSAEL (doted curve) and NA (dot-dashed
curves) when $n=100$.}\label{Fig:1}
\label{fig1}\end{figure}

The confidence regions of $\beta_0$ and their coverage
probabilities, with nominal level $1-\alpha=0.95$, were computed
from 500 runs. Four methods were used to construct the confidence
regions: the estimated empirical likelihood (EEL) with a
conditional approximation, the adjusted empirical likelihood
(AEL), the improved Rao--Scott adjusted empirical likelihood (IRSAEL)
and the normal approximation (NA). A comparison among three
methods was made through coverage accuracies and coverage areas of
the confidence regions. The simulated results are given in Figure~\ref{fig1}.

From Figure~\ref{fig1} we can see that EEL, AEL and IRSAEL give smaller confidence
regions than NA, and the region obtained by AEL is much smaller
than the others. Thus, AEL is the best of the four algorithms.

\begin{figure}
\begin{tabular}{cc}

\includegraphics{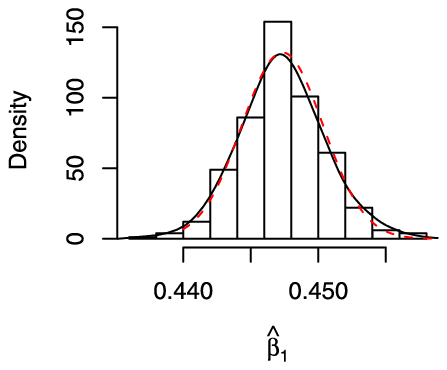}
&\includegraphics{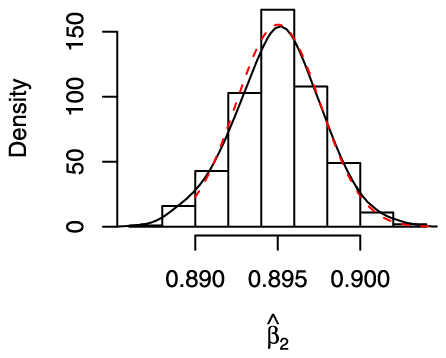}\\
(a) Histogram&(b) Histogram\\[6pt]

\includegraphics{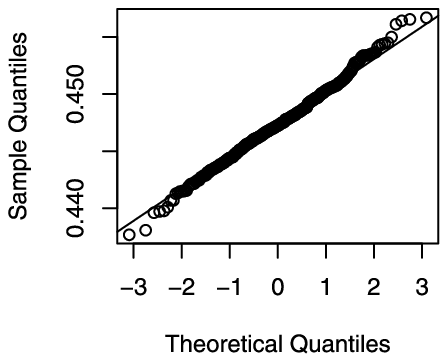}
&\includegraphics{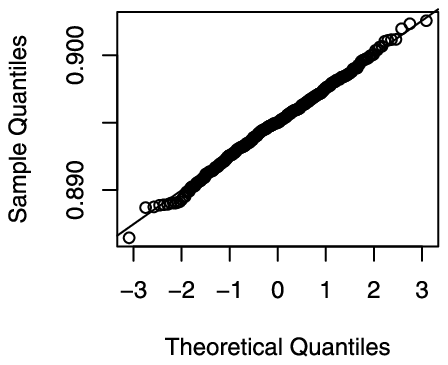}\\
(c) Normal Q--Q Plot&(d) Normal Q--Q Plot
\end{tabular}
\caption{(a) for $\beta_1$ and
(b) for $\beta_2$: the histograms of the 500 estimators of every
parameter, the estimated curve of density (solid curve) and the
curve of normal density (dashed curve); (c)~for $\beta_1$ and (d)
for $\beta_2$: the Q--Q plot of the 500 estimators of every
parameter.}\label{Fig:2}
\label{fig2}\end{figure}

The histograms of the 500 estimators of the parameter $\beta_1$
and $\beta_2$ are in Figures~\ref{fig2}(a) and (b), respectively. The Q--Q
plots of the 500 estimators of the parameter $\beta_1$ and
$\beta_2$ are in Figures~\ref{fig2}(c) and (d), respectively.

Figure~\ref{fig2} shows empirically that these estimators are
asymptotically normal. The means of the estimates of the unknown
parameters $\beta_1$ and $\beta_2$ are 0.44734 and 0.89502,
respectively, and their biases (standard deviations) are 0.000131
(0.00302) and
0.000596 (0.00257), respectively.

\begin{figure}
\begin{tabular}{cc}

\includegraphics{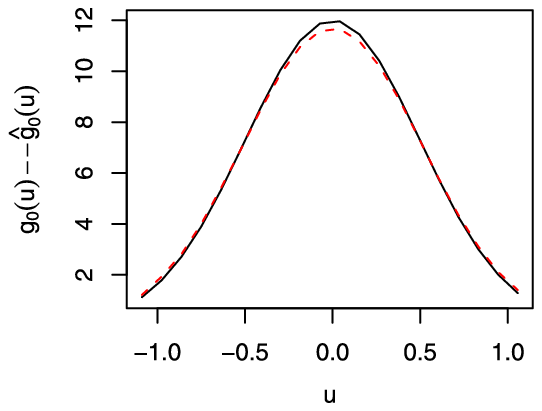}
&\includegraphics{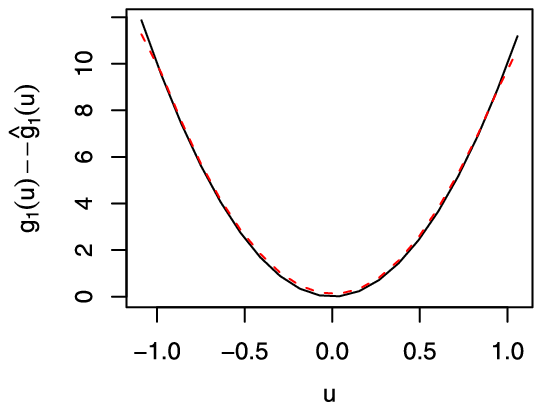}\\
(a)&(b)\\[5pt]

\includegraphics{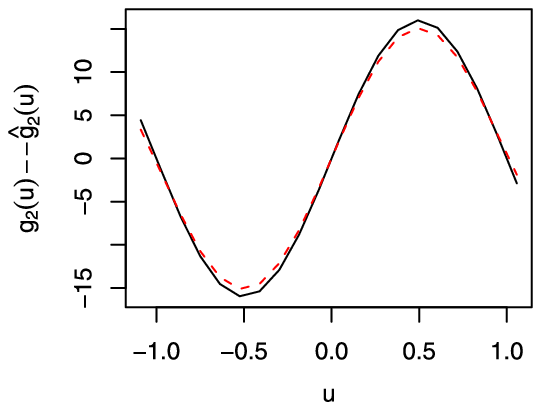}
&\includegraphics{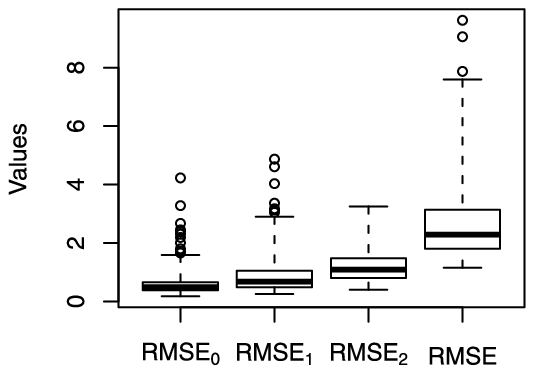}\\
(c)&(d)
\end{tabular}
 \caption{The true cure (solid
curve) and the estimated curve (dashed curve). (a) for
$g_{0}(\cdot)$, (b) for $g_{1}(\cdot)$, (c) for $g_{2}(\cdot)$; (d)
the boxplots of the 500 RMSE values in estimations of $g_0(\cdot)$,
$g_1(\cdot)$, $g_2(\cdot)$ and the sum of the three RMSEs.}
\label{Fig:3}
\label{fig3}\vspace*{-3pt}
\end{figure}

We also consider the average estimates of the coefficient
functions ${g}_{0}(u)$, ${g}_{1}(u)$ and ${g}_{2}(u)$ over the 500
replicates. The estimators $\hat{g}_j(\cdot)$ are assessed via the
root mean squared errors (RMSE); that is, $\mathrm{RMSE} = \sum_{j=0}^2{\rm
RMSE}_j$, where
\[
{\rm RMSE}_j= \Biggl[n_{\rm grid}^{-1}\sum_{k=1}^{n_{\rm grid}}
\{\hat{g}_j(u_k)-g_j(u_k)\}^2 \Biggr]^{1/2},
\]
and $\{u_k, k=1,\ldots,n_{\rm grid}\}$ are regular grid points.
The boxplot for the 500 RMSEs is given in Figure~\ref{fig3}.

From Figures~\ref{fig3}(a)--(c) we see every estimated curve agrees with
the true function curve very closely. Figure~\ref{fig3}(d) shows that all
RMSEs of estimates for the unknown functions are very small.\vspace*{-3pt}


\begin{appendix}\label{appm}
\section*{Appendices}\vspace*{-3pt}

We divide the appendices into Appendix~\ref{appmA} and
Appendix~\ref{appmB}. The proofs of \mbox{Theorems~\ref{theo1}--\ref{theo4}} are presented in Appendix~\ref{appmA},
and the proofs of Lemmas~\ref{lemm2} and~\ref{lemm3} are presented in Appendix~\ref{appmB}. We
use $c$ to represent any positive constant which may take a
different value for each appearance.\vspace*{-3pt}

\setcounter{section}{0}
\section{Proofs of theorems}\label{appmA}
\vspace*{-3pt}

The following lemma gives uniformly convergent rates of
$\hat{{\rm g}}(u;\beta)$ and $\hat{\dot{{\rm g}}}(u;\beta)$. This
lemma is
a straightforward extension of known results in nonparametric
function estimation; for its proof, the reader may refer to
Theorem~\ref{theo2}
in Wang and Xue \cite{WanXue11}, we hence omit the proof.\vadjust{\goodbreak}

\begin{lemma}\label{lemm1}  Suppose that conditions {\rm(C1)--(C3)},
{\rm(C5)} and {\rm(C6)} hold. Then
\[
\sup_{u\in{\cal U}_w,\beta\in{\cal
B}_n}\|\hat{{\rm g}}(u;\beta)-{{\rm g}}_0(u)\| =
\mathrm{O}_P \biggl(\biggl \{\frac{\log(1/h)}{nh} \biggr\}^{1/2} + h^2  \biggr)
\]
 and\vspace*{-2pt}
\[
\sup_{u\in{\cal U}_w,\beta\in{\cal
B}_n}\|\hat{\dot{{\rm g}}}(u;\beta)-\dot{{\rm g}}_0(u)\| =
\mathrm{O}_P \biggl(\biggl \{\frac{\log(1/h)}{nh^3} \biggr\}^{1/2} + h  \biggr).
\]
\end{lemma}

 Denote ${\cal G}=\{{{\rm g}}\dvt {\cal U}_w\times{\cal
B}\mapsto R^q\}$, $\|{{\rm g}}\|_{\cal G}=\sup_{u\in{\cal
U}_w,\beta\in{\cal B}_n}\|{{\rm g}}(u;\beta)\|$. From Lemma~\ref{lemm1}, we
have $\|\hat{{\rm g}}-{{\rm g}}_0\|_{\cal G}=\mathrm{o}_P(1)$ and
$\|\hat{\dot{{\rm g}}}-\dot{{\rm g}}_0\|_{\cal G}=\mathrm{o}_P(1)$; hence we can
assume that ${{\rm g}}$ lies in ${\cal G}_{\delta}$ with
$\delta=\delta_n\rightarrow0$ and $\delta>0$, where\vspace*{-1pt}
%
\begin{equation}
{\cal G}_{\delta}=\{{{\rm g}}\in{\cal G}\dvt \|{{\rm g}}-{{\rm g}}_0\|
_{\cal G}\leq
\delta,\|\dot{{\rm g}}-\dot{{\rm g}}_0\|_{\cal G}\leq\delta\}.
\label{eqA.1}
\end{equation}
Let
${\rm g}_0(\beta^TX;\beta)=E\{{\rm g}_0(\beta_0^TX)|\beta^TX\}$ and $\dot{\rm g}_0(\beta^TX;\beta)=E\{\dot{\rm g}_0(\beta_0^TX)|\beta^TX\}$,
%
\begin{eqnarray}
\label{eqA.2}
Q({{\rm g}},\beta) & =&E[\{Y-{{\rm g}}^T(\beta^TX;\beta)Z\}\dot{{\rm
g}}{}^T(\beta^TX;\beta)ZXw(\beta^TX)],
\\
\label{eqA.3}
Q_n({{\rm g}},\beta) & =&\frac{1}{n}\sum_{i=1}^n\{Y_i-{{\rm
g}}^T(\beta^TX_i;\beta)Z_i\}\dot{{\rm g}}{}^T(\beta^TX_i;\beta
)Z_iX_iw(\beta^TX_i).
\end{eqnarray}

The following two lemmas are required for obtaining the proofs
of the theorems; their proofs can be found in Appendix~\ref{appmB}.

\begin{lemma}\label{lemm2}  Suppose that conditions {\rm(C1)--(C6)} hold.
Then
%
\begin{eqnarray}
\label{eqA.4}
  \sup_{({{\rm g}},\beta)\in{\cal G}_{\delta}\times{\cal B}_n}\|
J_1({{\rm g}},\beta)\| &=& \mathrm{o}_P(n^{-1/2}),
\\
\label{eqA.5}
  \sup_{\beta\in{\cal B}_n} \|J_2(\hat{{\rm g}},\beta)\| &=& \mathrm{o}_P(n^{-1/2}),
\\
\label{eqA.6}
  \sup_{({{\rm g}},\beta)\in{\cal G}_{\delta}\times{\cal B}_n}\|
J_3({{\rm g}},\beta)\| &=& \mathrm{o}(n^{-1/2}),
\\
\label{eqA.7}
  \sqrt{n}J_4(\hat{{\rm g}},\beta_0)&\stackrel{D}{\longrightarrow
}&N(0,\sigma^2{\A}(\beta_0)),
\end{eqnarray}
  where ${\A}(\beta_0)$ is defined in (\ref{eq2.6}),
%
\begin{eqnarray*}
J_1({{\rm g}},\beta)
& =&Q_n({{\rm g}},\beta)-Q({{\rm g}},\beta)-Q_n({{\rm g}}_0,\beta_0),
\\
J_2({{\rm g}},\beta)
& =& Q({{\rm g}},\beta)-Q({{\rm g}}_0,\beta) \\
&&{}   -\varpi({{\rm g}}_0(\beta^TX; \beta);\beta)\{{{\rm g}}(\beta
^TX;\beta)-{{\rm g}}_0(\beta^TX; \beta)\},
\\
J_3({{\rm g}},\beta)
& =&\varpi({{\rm g}}_0(\beta^TX),\beta)\{{{\rm g}}(\beta^TX;\beta
)-{{\rm g}}_0(\beta^TX)\} \\
&&{} -\varpi({{\rm g}}_0(\beta_0^TX; \beta),\beta_0)\{{{\rm g}}(\beta
_0^TX;\beta_0)-{{\rm g}}_0(\beta_0^TX; \beta)\}
\end{eqnarray*}
  and
\[
J_4(\beta_0,{{\rm g}})=Q_n({{\rm g}}_0,\beta_0) +
\varpi({{\rm g}}_0(\beta_0^TX),\beta_0)\{{{\rm g}}(\beta_0^TX;\beta
_0)-{{\rm g}}_0(\beta_0^TX)\}.
\]
\end{lemma}

\begin{lemma}\label{lemm3}  Suppose that conditions {\rm(C1)--(C6)} hold.
Then
%
\begin{eqnarray}
\label{eqA.8}
  \sup_{\beta\in{\cal B}_n}\|Q_n(\hat{{\rm g}},\beta)\| &=& \mathrm{O}_P(n^{-1/2}),
\\
\label{eqA.9}
  \sup_{\beta\in{\cal B}_n}\|{\R}_n(\beta)- \sigma^2{\B}(\beta
_0)\|&=& \mathrm{o}_P(1),
\\
\label{eqA.10}
  \sup_{\beta\in{\cal B}_n}\max_{1\leq i\leq n}\|\hat{\eta
}_i(\beta)\| &=& \mathrm{o}_P(n^{1/2}),
\\
\label{eqA.11}
  \sup_{\beta\in{\cal B}_n}\|\lambda(\beta)\| &=& \mathrm{o}_P(n^{-1/2}),
\end{eqnarray}
where $Q_n(\hat{{\rm g}},\beta)$ is defined in {\rm(\ref{eqA.3})},
${\R}_n(\beta)=n^{-1}\sum_{i=1}^n\hat{\eta}_i(\beta)\hat{\eta
}_i^T(\beta)$,
${\B}(\beta_0)$ is defined in condition {\rm(C7)} and
$\hat{\eta}_i(\beta)$ is defined in (\ref{eq2.2}).
\end{lemma}

\begin{pf*}{Proof of Theorem~\ref{theo1}}
Note that, when $\beta=\beta_0$, Lemma~\ref{lemm3}
also holds. Applying the Taylor expansion to (\ref{eq2.2}) and invoking
Lemma~\ref{lemm3}, we can obtain
%
\begin{equation}
-2\log\hat{L}(\beta_0) =
-\sum_{i=1}^n \biggl[\lambda^T\hat{\eta}_i(\beta_0)-\frac
{1}{2} \{\lambda^T\hat{\eta}_i(\beta_0) \}^2 \biggr]+\mathrm{o}_P(1).
\label{eqA.12}
\end{equation}
%
By (\ref{eq2.3}) and Lemma~\ref{lemm3}, we have
\[
\sum_{i=1}^n \{\lambda^T\hat{\eta}_i(\beta_0)\}^2 = \sum
_{i=1}^n\lambda^T\hat{\eta}_i(\beta_0)+\mathrm{o}_P(1)
\]
and
\[
\lambda= \Biggl\{\sum_{i=1}^n\hat{\eta}_i(\beta_0)\hat{\eta
}_i^T(\beta_0)
 \Biggr\}^{-1} \sum_{i=1}^n \hat{\eta}_i(\beta_0) +
\mathrm{o}_P (n^{-1/2} ).
\]
This together with (\ref{eqA.12}) proves that
%
\begin{equation}
-2\log\hat{L}(\beta_0) =
nQ_n^T(\hat{{\rm g}},\beta_0){\R}_n^{-1}(\beta_0)Q_n(\hat{{\rm
g}},\beta_0)+\mathrm{o}_P(1),
\label{eqA.13}
\end{equation}
where $Q_n(\hat{{\rm g}},\beta_0)$ and ${\R}_n(\beta_0)$ are
defined in
(\ref{eqA.3}) and (\ref{eqA.9}), respectively. From (\ref{eqA.9}) of Lemma~\ref{lemm3} and (\ref{eqA.13}), we
obtain
%
\begin{equation}
\hspace*{-10pt}-2\log\hat{L}(\beta_0)\!=\!
 \bigl\{(\sigma^2{\A})^{-\fraca{1}{2}}\sqrt{n}Q_n(\hat{{\rm
g}},\beta_0) \bigr\}^T{\G}(\beta_0)
 \bigl\{(\sigma^2{\A})^{-\fraca{1}{2}}\sqrt{n}Q_n(\hat{{\rm
g}},\beta_0) \bigr\}\!+\!\mathrm{o}_P(1),
\label{eqA.14}
\end{equation}
where
${\G}(\beta_0)={\A}^{1/2}(\beta_0){\B}^{-1}(\beta_0){\A
}^{1/2}(\beta_0)$.
Let ${\G}_0 = \operatorname{diag}(w_1,\ldots,w_p)$, where $w_i$, \mbox{$1\leq
i\leq p$}, are the eigenvalues of ${\G}(\beta_0)$. Then there
exists an orthogonal matrix ${\Hbf}$ such that
${\Hbf}^T{\G}_0{\Hbf}={\G}(\beta_0)$. Using the notations of Lemma~\ref{lemm2}, we have
%
\begin{equation}
Q_n(\hat{{\rm g}},\beta) = J_1(\hat{{\rm g}},\beta) + J_2(\hat
{{\rm g}},\beta) +
J_3(\hat{{\rm g}},\beta) + J_4(\hat{{\rm g}},\beta_0) + Q({{\rm
g}}_0,\beta).
\label{eqA.15}
\end{equation}
Noting that $Q({{\rm g}}_0,\beta_0)=0$, from the above equation and
Lemma~\ref{lemm2}, we have
\[
Q_n(\hat{{\rm g}},\beta_0)=J_4(\hat{{\rm g}},\beta_0) + \mathrm{o}_P(n^{-1/2}).
\]
Hence, by (\ref{eqA.7}) of Lemma~\ref{lemm2}, we have
\[
{\Hbf}\{\sigma^{-2}{\A}^{-}(\beta_0)\}^{1/2}\sqrt{n}Q_n(\hat{{\rm
g}},\beta_0)
\stackrel{{D}}{\longrightarrow} N(0,{\I}_{p}),
\]
where ${\I}_{p}$ is the $p\times p$ identity matrix. This together
with (\ref{eqA.14}) proves Theorem~\ref{theo1}.
\end{pf*}

\begin{pf*}{Proof of Theorem~\ref{theo2}} Under the conditions of Theorem~\ref{theo2}, we
can follow similar arguments to those used by   Wang and Xue \cite{WanXue11} and show that $\hat
{\beta}$ is
a root-$n$ consistent estimator of $\beta_0$. Because the proof is
straightforward, we do not present it here. We next demonstrate
the asymptotic normality of $\hat{\beta}$. By Lemma~\ref{lemm3} and, similarly
to the proof of~(\ref{eqA.13}), we can obtain
%
\begin{equation}
\log\hat{L}(\beta) =
-\frac{n}{2}Q_n^T(\hat{{\rm g}},\beta)\{\sigma^2{\B}(\beta)\}
^{-1}Q_n(\hat{{\rm g}},\beta)+\mathrm{o}_P(1),
\label{eqA.16}
\end{equation}
uniformly for $\beta\in{\cal B}_n$, where $\mathrm{o}_P(1)$ tends to 0 in
probability uniformly for $\beta\in{\cal B}_n$. Since the
estimator $\hat\beta$ is a maximum of $\log\hat{L}(\beta)$, and
${\B}(\beta_0)$ is a positive definite matrix, the resulting
estimator $\hat{\beta}$ is equivalent to solving the estimation
equation $Q_n(\hat{{\rm g}},\beta)=0$; that is,
$Q_n(\hat{{\rm g}},\hat{\beta})=0$. Note that $Q({{\rm g}}_0,\beta
_0)=0$, and we
then have, by Taylor's expansion, that
%
\begin{equation}
Q({{\rm g}}_0,\beta) = -{\B_*}(\beta_0)(\beta-\beta_0)+\mathrm{o}(n^{-1/2}),
\label{eqA.17}
\end{equation}
uniformly for $\beta\in{\cal B}_n$, where ${\B_*}(\beta_0)$ is the
same as that in (\ref{eqA.9}). By (\ref{eqA.15}), (\ref{eqA.17}) and~(\ref{eqA.4})--(\ref{eqA.6}) of Lemma~\ref{lemm2},
we have
\[
Q_n(\hat{{\rm g}},\hat{\beta})=J_4(\hat{{\rm g}},\beta_0) -
{\B_*}(\beta_0)(\hat{\beta}-\beta_0) + \mathrm{o}_P(n^{-1/2}).
\]
Noting that $Q_n(\hat{{\rm g}},\hat{\beta})=0$, we get
\[
\sqrt{n} (\hat{\beta}-\beta_0 ) =
\sqrt{n} \B_*^{-1}(\beta_0)J_4(\hat{{\rm g}},\beta_0) + \mathrm{o}_P(1).
\]
This together with (\ref{eqA.7}) of Lemma~\ref{lemm2} proves Theorem~\ref{theo2}.
\end{pf*}

\begin{pf*}{Proof of Theorem~\ref{theo3}}   Decomposing $\hat{\sigma}^2$ into
several parts, we get
%
\begin{eqnarray*}
\hat{\sigma}^2
& =& \frac{1}{n}\sum_{i=1}^n\varepsilon_i^2
+ \frac{1}{n}\sum_{i=1}^n [{{\rm g}}_0(X_i^T\beta_0)-\hat{{\rm
g}}(X_i^T\hat{\beta};\hat{\beta})\}^TZ_i ]^2\\
&&{} + \frac{2}{n}\sum_{i=1}^n\varepsilon_i \{{{\rm
g}}_0(X_i^T\beta_0)-\hat{{\rm g}}(X_i^T\hat{\beta};\hat{\beta
}) \}^T Z_i\\
& \equiv& I_1 + I_2 + I_3.
\end{eqnarray*}
Using the central limit theorem, we have
\[
\sqrt{n}(I_1-\sigma^2)=\frac{1}{\sqrt{n}}\sum_{i=1}^n(\varepsilon
_i^2-\sigma^2)
\stackrel{D}{\longrightarrow}N(0,\operatorname{var}(\varepsilon^2)).
\]
By Lemma~\ref{lemm1}, we can obtain
\[
|I_2|
\leq\frac{1}{n}\sum_{i=1}^n\|Z_i\|^2
 \Bigl\{\sup_{(u,\beta)\in({\cal U}_w,{\cal B}_n)}\|\hat{{\rm
g}}(u;\beta)-{{\rm g}}_0(u)\| \Bigr\}^2
= \mathrm{o}_P (n^{-1/2} ).
\]
For $I_3$, we have
%
\begin{eqnarray*}
I_3 & =& \frac{2}{n}\sum_{i=1}^n\varepsilon_i \{{{\rm
g}}_0(X_i^T\beta_0)-\hat{{\rm g}}(X_i^T\beta_0;\beta_0) \}^T
Z_i\\
& &{} + \frac{2}{n}\sum_{i=1}^n\varepsilon_i \{\hat{{\rm
g}}(X_i^T\beta_0;\beta_0)-\hat{{\rm g}}(X_i^T\hat{\beta};\hat
{\beta}) \}^T Z_i\\
& \equiv& I_{31} + I_{32}.
\end{eqnarray*}
It is not hard to show that $I_{31}=\mathrm{O}_P (n^{-1/2} )$. By
Theorems~\ref{theo1} and~\ref{theo3}, we obtain
\[
 |I_{32}|\leq\frac{2}{n}\sum_{i=1}^n\bigl(\|Z_i\||\varepsilon_i|\|X_i-E(X_i|\beta_0^TX_i)\|\bigr)\|\hat{\beta}-\beta_0\|\mathrm{O}_P(1)
    = \mathrm{O}_P(n^{-1/2}).
\]
This together with above results proves Theorem~\ref{theo3}.
\end{pf*}

\begin{pf*}{Proof of Theorem~\ref{theo4}}  Note that
$\hat{\A}(\beta_0)\stackrel{P}{\longrightarrow}{\A}(\beta_0)$ and
$\hat{\B}(\beta_0)\stackrel{P}{\longrightarrow}{\B}(\beta_0)$. By
the expansion of $\hat{l}_{\rm ael}(\beta_0)$, defined in (\ref{eq3.1}) and
(\ref{eqA.16}), we get
%
\begin{equation}
\hat{l}_{\rm ael}(\beta_0) =
nQ_n^T(\hat{{\rm g}},\beta_0)\{\sigma^{-2}{\A}^{-}(\beta_0)\}
Q_n(\hat{{\rm g}},\beta_0)+\mathrm{o}_P(1).
\label{eqA.18}
\end{equation}
This together with (\ref{eqA.15}) and (\ref{eqA.18}) proves Theorem~\ref{theo4}.
\end{pf*}

\section{Proofs of lemmas}\label{appmB}

\begin{pf*}{Proof of Lemma~\ref{lemm2}}  We first prove
(\ref{eqA.4}). Denote
$r_n({{\rm g}},\beta)=\sqrt{n}\{Q_n({{\rm g}},\beta)-Q({{\rm
g}},\beta)\}$.
Noting that $Q({{\rm g}}_0,\beta_0)=0$, we clearly have
%
\begin{equation}
J_1({{\rm g}},\beta) = n^{-1/2}\{r_n({{\rm g}},\beta)-r_n({{\rm
g}}_0,\beta_0)\}.
\label{eqB.1}
\end{equation}
It can be shown that the empirical process
$\{r_n({{\rm g}},\beta)\dvt {{\rm g}}\in{\cal G}_1,\beta\in{\cal B}_1\}
$ has the
stochastic equicontinuity, where ${\cal B}_1=\{\beta\in{\cal B}\dvt
\|\beta-\beta_0\|\leq1\}$ and ${\cal G}_1$ are defined in (\ref{eqA.1})
with $\delta=1$, which are subsets of ${\cal B}$ and ${\cal G}$,
respectively. The equicontinuity is sufficient for proof of (\ref{eqA.4})
since $\delta<1$ for large enough $n$. This stochastic
equicontinuity follows by checking the conditions of Theorem~\ref{theo1} in
Doukhan, Massart and Rio \cite{DouMasRio95}. Therefore, we have
$r_n({{\rm g}},\beta)-r_n({{\rm g}}_0,\beta_0)=\mathrm{o}_P(1)$, uniformly for
$\beta\in{\cal B}_1$ and ${{\rm g}}\in{\cal G}_1$. This together with~(\ref{eqB.1}) proves (\ref{eqA.4}).

We now prove (\ref{eqA.5}). Define the functional derivative $\varpi({\rm g}_0(\cdot;\beta),\beta)$ of
$Q({{\rm g}},\beta)$ with respect to ${{\rm g}}(\cdot;\beta)$ at
${\rm g}_0(\cdot;\beta)$ at the direction ${{\rm g}}(\cdot;\beta)-{\rm g}_0(\cdot;\beta)$
by
%
\begin{eqnarray*}
&& \varpi({\rm g}_0(\cdot;\beta),\beta)\{{{\rm g}}(\cdot;\beta)-{\rm g}_0(\cdot;\beta)\} \\[-1pt]
&& \quad   = \lim_{\tau\rightarrow
0}\bigl[Q\bigl({\rm g}_0(\cdot;\beta)+\tau\bigl({{\rm g}}(\cdot;\beta)-{\rm g}_0(\cdot;\beta)\bigr),\beta\bigr) -
Q({\rm g}_0(\cdot;\beta),\beta)\bigr]\cdot\frac{1}{\tau},
\end{eqnarray*}
where $Q({{\rm g}},\beta)$ is defined in (\ref{eqA.2}). We have
%
\begin{eqnarray}
\label{eqB.2}
&& \varpi({{\rm g}}_0(\beta^TX; \beta),\beta)\{{{\rm g}}(\beta^TX;\beta
)-{{\rm g}}_0(\beta^TX; \beta)\} \nonumber
\\[-8.5pt]
\\[-8.5pt]
&& \quad     = - E[\{{{\rm g}}(\beta^TX;\beta)-{{\rm g}}_0(\beta^TX; \beta)\}
^TZ\dot{{\rm g}}_0^T(\beta^TX; \beta)ZXw(\beta^TX)].\nonumber
\end{eqnarray}
It follows from (\ref{eqB.2}) that
%
\begin{eqnarray*}
J_2({{\rm g}},\beta)
& =& - E [\{{{\rm g}}(\beta^TX;\beta)-{{\rm g}}_0(\beta_0^TX)\}
^TZXZ^T \\[-1pt]
&&\hphantom{- E [} {}
\times\{{\dot{{\rm g}}}(\beta^TX;\beta)-{\dot{{\rm g}}}_0(\beta
^TX; \beta)\}w(\beta^TX) ],
\end{eqnarray*}
and hence we have
%
\begin{eqnarray}
\label{eqB.3}
\omega^TJ_2(\hat{{\rm g}},\beta)
& =& - \int\{\hat{{\rm g}}(u;\beta)-{{\rm g}}_0(u)\}^T\mu_\omega(u)
\nonumber
\\[-8.5pt]
\\[-8.5pt]
& &\hphantom{- \int}{} \times\{\hat{\dot{{\rm g}}}(u;\beta)-{\dot{{\rm
g}}}_0(u)\}w(u)f(u)\,\mathrm{d}u
+\mathrm{o}_P(n^{-1/2})\nonumber
\end{eqnarray}
for any $p$-dimension vector $\omega$, where
$\mu_\omega(u)=E\{Z\omega^TXZ^T|\beta^TX=u\}$, and $f(u)$ is the
probability density of $\beta^TX$. Using the standard argument of
nonparametric estimation, we can prove
%
\begin{equation}
\hat{{\rm g}}(u;\beta)-{{\rm g}}_0(u)
= {\D}^{-1}(u)\{f(u)\}^{-1}\xi_{n}(u;\beta)
+ \mathrm{O}_P(c_n),
\label{eqB.4}
\end{equation}
uniformly for $u\in{\cal U}_w$ and $\beta\in{\cal B}_n$, where
$c_n=n^{-1/2}+h^2$ and ${\D}(u)$ is defined in condition~(C6).
\[
\xi_{n}(u;\beta) =
\frac{1}{n}\sum_{i=1}^nZ_i\{Y_i-{{\rm g}}_0^T(\beta^TX_i)Z_i\}
K_h(\beta^TX_i-u).
\]
This together with (\ref{eqB.3}) derives that
%
\begin{eqnarray*}
\omega^TJ_2(\hat{{\rm g}},\beta)
& =& - \int\{{\D}^{-1}(u)\xi_{n}(u;\beta)\}^T
\mu_\omega(u)\{\hat{\dot{{\rm g}}}(u;\beta)-{\dot{{\rm g}}}_0(u)\}\,\mathrm{d}u + \mathrm{O}_P(c_n) \\[-1pt]
& =& -n^{-1/2}\{\gamma_n(\hat{\dot{{\rm g}}},\beta)-\gamma_n(\dot
{{\rm g}}_0,\beta)\}+ \mathrm{O}_P(c_n),
\end{eqnarray*}
where
$\gamma_n({\dot{{\rm g}}},\beta)=n^{-1/2}\sum_{i=1}^n\varepsilon
_iw(\beta^TX_i)Z_i^T{\D}^{-1}(\beta^TX_i)\mu_\omega(\beta
^TX_i){\dot{{\rm g}}}(\beta^TX_i;\beta)$.
Using the empirical process techniques, and similarly to the proof
of (\ref{eqA.4}), we can show that the stochastic equicontinuity of
$\gamma_n({\dot{{\rm g}}},\beta)$, and hence
$\|\gamma_n(\hat{\dot{{\rm g}}},\beta)-\gamma_n(\dot{{\rm
g}}_0,\beta)\|=\mathrm{o}_P(1)$.
Also, $nh^4=\mathrm{O}(1)$ implies $h^2=\mathrm{O}(n^{-1/2})$, and hence
$c_n=\mathrm{O}(n^{-1/2})$. Thus, the proof of (\ref{eqA.5}) is complete.

We now prove (\ref{eqA.6}). Denote
$\psi({\dot{{\rm g}}_0},\beta)=\dot{{\rm g}}_0^T(\beta^TX;\beta)ZXw(\beta^TX)$ and
$\varphi({{\rm g}},\beta)=\{{{\rm g}}(\beta^TX; \beta)-{{\rm
g}}_0(\beta^TX; \beta)\}^TZ$. It
follows from (\ref{eqB.2}) that
%
\begin{eqnarray*}
J_3({{\rm g}},\beta)
& =& - E\{\varphi({{\rm g}},\beta)\psi({\dot{{\rm g}}_0},\beta)\} + E\{
\varphi({{\rm g}},\beta_0)\psi({\dot{{\rm g}}_0},\beta_0)\} \\
& =&
- E[\{\varphi({{\rm g}},\beta)-\varphi({{\rm g}},\beta_0)\}\psi
({\dot{{\rm g}}_0},\beta)] \\
&&{} - E[\varphi({{\rm g}},\beta_0)\{\psi({\dot{{\rm g}}_0},\beta
)-\psi({\dot{{\rm g}}_0},\beta_0)\}] \\
& \equiv& J_{31}({{\rm g}},\beta) + J_{32}({{\rm g}},\beta).
\end{eqnarray*}
By condition (C2), we get
%
\begin{eqnarray*}
&& \|\varphi({{\rm g}},\beta)-\varphi({{\rm g}},\beta_0)\|
\\
 && \quad = \|[\{{\rm g}(\beta^TX;\beta)-{\rm g}(\beta_0^TX;\beta_0)\}-\{{\rm g}_0(\beta^TX;\beta)-{\rm g}_0(\beta_0^TX)\}]^TZ\|
 \\
 && \quad = \|[\{\dot{\rm g}(\beta_1^TX;\beta_1)-\dot{\rm g}_0(\beta_2^TX)\}(\beta-\beta_0)^T\{X-E(X|\beta_0^TX)\}]^TZ\|
 \\
 && \quad \leq c\|\dot{\rm g}-\dot{\rm g}_0\|_{\cal G}\|\beta-\beta_0\|(\|X-E(X|\beta_0^TX)\|)(\|Z\|),
\end{eqnarray*}
where $\beta_1$ and $\beta_2$ are between $\beta$ and $\beta_0$,
and $\|\psi({\dot{\rm g}_0},\beta)\|\leq c(\|Z\|)(\|X\|)$. Therefore, we have
$\|J_{31}({{\rm g}},\beta)\| = \mathrm{o}(n^{-1/2})$, uniformly for
${{\rm g}}\in{\cal G}_{\delta}$ and $\beta\in{\cal B}_n$.
Similarly, we
can prove $\|J_{32}({{\rm g}},\beta)\|=\mathrm{o}(n^{-1/2})$, uniformly for
${{\rm g}}\in{\cal G}_{\delta}$ and $\beta\in{\cal B}_n$, and hence
(\ref{eqA.6}) follows.

Finally, we prove (\ref{eqA.7}). Let $f_0(u)$ denote the density function
of $\beta_0^TX$. By (\ref{eqB.2}) and~(\ref{eqB.4}), and using the dominated
convergence theorem (Lo\`{e}ve \cite{Loe00}), we can obtain
%
\begin{eqnarray*}
&& \varpi({{\rm g}}_0(\beta_0^TX),\beta_0)\{\hat{{\rm g}}(\beta
_0^TX;\beta_0)-{{\rm g}}_0(\beta_0^TX)\}
\\
&& \quad   = - \int{\C}(u)\{\hat{{\rm g}}(u;\beta_0)-{{\rm g}}_0(u)\}f_0(u)\,\mathrm{d}u
\\
&& \quad   = -
\frac{1}{n}\sum_{i=1}^n\varepsilon_i{\C}(\beta_0^TX_i){\D
}^{-1}(\beta_0^TX_i)Z_i  + \mathrm{o}_P(c_n).
\end{eqnarray*}
This together with (\ref{eqA.3}) proves that
\[
J_4(\hat{{\rm g}},\beta_0) =
\frac{1}{n}\sum_{i=1}^n\varepsilon_i\zeta_i + \mathrm{o}_P(c_n),
\]
where $\zeta_i=V_i - {\C}(\beta_0^TX_i){\D}^{-1}(\beta_0^TX_i)Z_i$
and $V_i=X_i\dot{{\rm g}}_0^T(\beta_0^TX_i)Z_iw(\beta_0^TX_i)$.
Therefore, by the central limit theorem and Slutsky's theorem, we
get
\[
\sqrt{n}J_4(\hat{{\rm g}},\beta_0)= \frac{1}{\sqrt{n}}\sum
_{i=1}^n\varepsilon_i\zeta_i+\mathrm{o}_P(1)
\stackrel{D}{\longrightarrow}N(0,\sigma^2{\A}(\beta_0)).
\]
This proves (\ref{eqA.7}). The proof of Lemma~\ref{lemm2} is
complete.
\end{pf*}

\begin{pf*}{Proof of Lemma~\ref{lemm3}} By (\ref{eqA.15}), (\ref{eqA.17}) and Lemma~\ref{lemm2}, we can
prove (\ref{eqA.8}). We now prove (\ref{eqA.9}). Let
%
\begin{eqnarray*}
R_{ni}(\beta)
& =&
\varepsilon_i\dot{{\rm g}}_0^T(\beta_0^TX_i)Z_iX_i\{w(\beta
^TX_i)-w(\beta_0^TX_i)\}
\\
&&{} + \varepsilon_i\{\hat{\dot{{\rm g}}}(\beta^TX_i;\beta
)-\dot{{\rm g}}_0(\beta_0^TX_i)\}^TZ_iX_iw(\beta^TX_i)
\\
&&{} + \{{{\rm g}}_0(\beta_0^TX_i)-\hat{{\rm g}}(\beta
^TX_i;\beta)\}^TZ_iZ_i^T\dot{{\rm g}}_0(\beta_0^TX_i)X_iw(\beta^TX_i)
\\
&&{} + \{{{\rm g}}_0(\beta_0^TX_i)-\hat{{\rm g}}(\beta
^TX_i;\beta)\}^TZ_iZ_i^T
\\
&&{} \times\{\hat{\dot{{\rm g}}}(\beta^TX_i;\beta)-\dot{{\rm
g}}_0^T(\beta_0^TX_i)\}X_iw(\beta^TX_i).
\end{eqnarray*}
Then we have $\hat{\eta}_i(\beta)=\eta_i(\beta_0) +
R_{ni}(\beta)$, where $\eta_i(\cdot)$ is defined in (\ref{eq2.1}), and
hence
%
\begin{eqnarray}
\label{eqB.5}
{\R}_n(\beta) & =&
\frac{1}{n}\sum_{i=1}^n\eta_i(\beta_0)\eta_i^T(\beta_0)
+ \frac{1}{n}\sum_{i=1}^nR_{ni}(\beta)R_{ni}^T(\beta)
\nonumber\\
&&{} + \frac{1}{n}\sum_{i=1}^n\eta_i(\beta_0)R_{ni}^T(\beta)
+ \frac{1}{n}\sum_{i=1}^nR_{ni}(\beta)\eta_i^T(\beta_0)
\\
& \equiv& M_1(\beta_0) + M_2(\beta) + M_3(\beta) + M_4(\beta).\nonumber
\end{eqnarray}
By the law of large numbers, we have
$M_1(\beta_0)\stackrel{P}{\longrightarrow}\sigma^2{\B}(\beta_0)$.
Therefore, to prove (\ref{eqA.9}), we only need to show that
$M_k(\beta)\stackrel{P}{\longrightarrow}0$ uniformly for $\beta$,
$k=2,3,4$.

Let $M_{2,st}(\beta)$ denote the $(s,t)$ element of $M_2(\beta)$,
and $R_{ni,s}(\beta)$ denote the $s$th component of
$R_{ni}(\beta)$. Then by the Cauchy--Schwarz inequality, we have
%
\begin{equation}
|M_{2,st}(\beta)|\leq
 \Biggl(\frac{1}{n}\sum_{i=1}^nR_{ni,s}^2(\beta) \Biggr)^{1/2}
\Biggl(\frac{1}{n}\sum_{i=1}^nR_{ni,t}^2(\beta) \Biggr)^{1/2}.
\label{eqB.6}
\end{equation}
%
It can be shown by a direct calculation that
%
\[
\frac{1}{n}\sum_{i=1}^nR_{ni,s}^2(\beta)\stackrel
{P}{\longrightarrow}0,
\]
uniformly for $\beta\in{\cal B}_n$. This together with (\ref{eqB.6})
proves that $M_2(\beta)\stackrel{P}{\longrightarrow}0$, uniformly
for $\beta\in{\cal B}_n$. Similarly, it can be shown that
$M_3(\beta)\stackrel{P}{\longrightarrow}0$ and
$M_4(\beta)\stackrel{P}{\longrightarrow}0$, uniformly for
$\beta\in{\cal B}_n$. This together with (\ref{eqB.5}) proves (\ref{eqA.9}).


Similarly to above proof, we can derive (\ref{eqA.10}). (\ref{eqA.11}) can be
shown by using (\ref{eqA.8})--(\ref{eqA.10}), and employing the same arguments
used in the proof of (2.14) in Owen \cite{Owe90}.
\end{pf*}
\end{appendix}

\section*{Acknowledgements}
We are grateful for the many detailed suggestions of the editor, the
associate editor and the referees, which led to significant
improvements of the paper.
Liugen Xue's research was supported by the National Natural Science
Foundation of China (10871013,11171012), the Beijing Natural Science Foundation
(1102008), the Beijing Municipal Education Commission Foundation (KM201110005029) and the PHR(IHLB).
Qihua Wang's research was supported by the National Science Fund for
Distinguished Young Scholars in China (10725106), National Natural
Science Foundation of China (10671198), the National Science Fund for
Creative Research Groups in China and a~grant from the Key Lab of
Random Complex Structure and Data Science, CAS and
the Key grant from Yunnan Province (2010CC003).

%

\printhistory

\end{document}